# On permutations derived from integer powers $x^n$

John S. McCaskill[1,2,*] and Peter R. Wills[3]

[1]BioMIP Research Group, Ruhr University Bochum, Germany

[2] European Center for Living Technology, Universitas Ca' Foscari, Venice, Italy

[3] Department of Physics, University of Auckland, New Zealand

## ABSTRACT

We present a general theorem characterizing the relationship between the prime base $p$ representations of non-negative integers $x$ and their positive integer powers, $x^n$. For any positive integer $l$, the theorem establishes the existence of bijective mappings (permutations) between all $p^l$ members $x$ of each non-zero residue class mod $p$ satisfying $x < p^{l+1}$. These mappings are obtained as the integer part of $x^p p^{-\alpha}$ for a particular positive integer $\alpha$, depending on $n$ and $p$, called the *coding shift*, for which an explicit formula is given. For relatively prime $n$ and $p$, $\alpha = 1$ and the result follows directly from properties of the multiplicative group of invertible elements modulo $p^{l+1}$. We extend our result for general $n$ also to identify the coding shift required to obtain such bijective mappings for members of the zero residue class mod $p$, demonstrating that such bijective mappings (or encodings) can be found for any finite $l$ and for all positive integers $x < p^{l+1}$.

## 1 Introduction

In this contribution, we state and prove a general theorem characterizing the relationship between the prime base $p$ representations of non-negative integers $x$ and their positive integer powers, $x^n$.

For any positive integer $l$ (the radix) and given prime $p$, let $X_{l,p} = \{x': 0 \leq x' < p^l\}$, be the set of $l$ digit numbers base $p$. We seek to characterise the relationship between the base $p$ representation of all numbers $x$ in $X_{l+1,p}$ having remainder $r$ modulo $p$ (which are relatively prime to $p$ if $r \neq 0$) and that of their $n^{\text{th}}$ powers, $x^n$.

In particular we shall identify bijective mappings (permutations) on $X_{l,p}$ induced by the integer mapping $x^n$. This power mapping maps the $l$ non-terminal digits of the base $p$ representation of $x$ onto $l$ consecutive digits (starting at the digit $\alpha(p,n)$ in the base $p$ representation of $x^n$). We call the relevant digits in $x^n$ the *coding region*. We shall determine the starting digit position $\alpha(p,n)$ (the *coding shift* or "shift" for short) below. Because these mappings are invertible, these digits provide a coded representation of the first $p^l$ numbers. The coding region in the base $p$ representation of a power $x^n$, together with the remainder (terminal digit) $r = x \bmod p$, in any prime base $p$, thus uniquely determines the remaining digits in the power and the original number in the range $0 < x < p^{l+1}$ that was raised to the power (i.e. the discrete $n$th root).

---

[*] Corresponding author: <johnsmccaskill@gmail.com>

Our investigation of this relationship was motivated by understanding the complexity of representations associated with multiplication, in particular in connection with random number generation and the inversion of binary operations. The proof of the result is particularly simple for the case of a relatively prime power $n$ with $(n,p) = 1$, following from Euler's theorem [1,2] and the properties of cosets used to prove Lagrange's theorem [3], as we shall see. In this case, the result is known to have applications in coding theory in the context of discrete roots of integers. In the more general case which we shall prove, the result appears novel and demonstrates that the power function induces non-trivial and complete recodings of arbitrarily large sets of consecutive numbers, without duplication or omissions. *En route* to the proof, we establish a new result about the divisibility of binomial coefficients involving prime powers.

## 2   Power Coding Theorem

For any $n \in \mathbb{N}^+$ (positive integers) and prime number $p$, the base, we may write $n = qp^k$ where the greatest common divisor of $q$ and $p$, $\gcd(q,p) = 1$. Here $q$ is a positive integer and $k \in \mathbb{N}$ is the exact power of $p$ in $n$. In particular, in the most common case that the chosen prime base $p$ is not a factor of $n$, then $k = 0$. Note that for any $x \in X_{l+1,p}$, we may write $x = px' + r$ with $x' \in X_{l,p}$ and $0 \leq r < p$. We can now state the theorem.

<u>Theorem 1 (Power Coding Theorem):</u>

> *For all $n, l \in \mathbb{N}^+$ and prime $p$, there exist, for each integer residue $r$ with $0 < r < p$, discrete bijective mappings (permutations) $C_r: X_{l,p} \to X_{l,p}$ defined by*
>
> $$C_r(x') = \left\lfloor \frac{(px'+r)^n}{p^{\alpha(n,p)}} \right\rfloor \bmod p^l \qquad (1)$$
>
> *for $x' \in X_{l,p} = \{x': 0 \leq x' < p^l\}$. The digit shift $\alpha(n,p) \in \mathbb{N}^+$ appearing in $C_r$ is*
>
> $$\alpha(n,p) = 1 + k + \delta_{p,2}(1 - \delta_{k,0}) \qquad (2)$$
>
> *where $k$ is the exponent of $p$ in the prime factorisation of $n$ ($n = qp^k, p \nmid q$).*

Here, as usual, $\lfloor y \rfloor$ denotes the integer part (floor) of $y$, i.e. the largest integer $\leq y$ and $\delta_{ij}$ is the Kronecker delta function (1 if $i = j$, 0 otherwise). Note that except for the binary case $p = 2$, $\alpha(n,p) = 1 + k_n$, and that if, as is most common, $p$ is not a factor of $n$, $\alpha(n,p) = 1$ (even if $p = 2$). In the remaining special case that $p = 2$ and $n$ is divisible by exactly $k_n > 0$ powers of 2, $\alpha(n,p) = 2 + k_n$.

Actually, given any prime $p$, any positive integer $x$ can be written uniquely in the form $x = p^j(px_j' + r_j)$ with integers $0 < r_j < p$ and $j \geq 0$. The following corollary extends the Power Coding Theorem to sets of numbers $x$ with $j > 0$.

Corollary 1:

> For all $n, l \in \mathbb{N}^+$ and prime $p$, and for any $x > 0$, $j \in \mathbb{N}$, a set of bijective mappings $C_{r_j}: X_{l,p} \to X_{l,p}$ exists, arising from the nth powers of numbers from sets of the form $X_{l,p,r} = \{x = p^j(px_j' + r_j): x' \in X_{l,p}\}$ with $0 < r_j < p$, where $x' \in X_{l,p} = \{x': 0 \leq x' < p^l\}$. They have the form

$$C_{r_j}(x_j') = \left\lfloor \frac{(px_j' + r_j)^n}{p^{\alpha'(n,p)}} \right\rfloor \bmod p^l \qquad (3)$$

> where

$$\alpha'(n,p) = nj + 1 + k + \delta_{p,2}(1 - \delta_{k,0}). \qquad (4)$$

We refer to the Corollary as the Extended Power Coding Theorem and the "coding region" is defined as the $l$ base $p$ digits of $x^n$ (coefficients of increasing powers of $p$) starting at digit $\alpha'(n,p) + 1$, (i.e. the coefficient of the power $p^{\alpha'(n,p)}$). Corollary 1 is a direct consequence of the above theorem, noting that if $j$ is the highest power of $p$ dividing $x$, we may write $x = x_j p^j$ with $x_j = px_j' + r_j$ relatively prime to $p$ and so $x^n = x_j^n p^{jn}$. Applying the Power Coding Theorem to the set of numbers $x_j$ gives the result directly, with the additional shift $nj$.

The bijective nature of the mappings immediately implies also that all the non terminal digits of any number $x$ in its base $p$ representation are uniquely determined by the corresponding coding region of its $n$th power $x^n$.

Corollary 2:

> For any positive integer $x$ raised to the power $n$, and for any integer $l \geq 1$, the $l$ digit coding region base prime $p$, defined by the Extended Power Coding Theorem, together *with the lowest non-zero digit, the residue $r \bmod p$, uniquely identifies the first $l + 1$ base $p$ digits of $x$, and so determines the integer nth root of $x^n$ up to multiples of $p^{l+1}$.*

Corollary 2 follows from the existence of an inverse mapping, guaranteed by the bijective nature of the mapping defined in the Extended Power Coding Theorem.

We provide some illustrations of the power coding theorem in the next section. The bijective mappings defined by the powers possess interesting structure as graphs for smaller powers base 2, and provide some remarkably simple pseudorandom permutations even for only slightly larger powers and bases. In section 4 we establish a preliminary result, which simplifies the proof of the Power Coding Theorem. While this result can be deduced from Legendre's [4] factorization of the factorials or more directly from Kummer's Theorem [5] on the exact power of $p$ dividing a general binomial coefficient, it appears to be less well known, but nevertheless quite useful, and so we provide an independent proof from first principles. In section 5, we then complete the proof of the Power Coding Theorem. In the final section we discuss some implications of this result.

## 3   Examples of the Power Coding Theorem

There is a straightforward graphical means to understand the Power Coding Theorem for the various relationships between the power $n$ and base $p$ as shown in fig. 1. In the diagrams, the boxes are placeholders for single digits in the base $p$ representation of the numbers $x$ and $z = x^n$. Boxes with …| as contents stand for multiple digits which are not depicted explicitly. In each case, the digits $z'$ base $p$ are labelled correctly only for the lowest exponent value $k$ of the base $p$ depicted. For higher exponents, the labelling of the coding region shifts to the left, as the shaded mapping and curly braces indicate. Not shown are the extended PCT bijective mappings, which are similar but shifted further to the left for sets of numbers divisible by powers of the base $p$.

As an example, consider the base $p = 3$ representations of all positive numbers $x < 3^3 = 27$. We consider separately the three sets of numbers divisible by $3^2$, 3 and not divisible by 3, and divide each of these sets further into the two non-zero residue classes 1 and 2 mod 3. These six sets partition all the numbers from 1 to 26.

$3 \nmid x$  $X_{3,2} = \{0,1,2,3,4,5,6,7,8\}$   $r = 1: x \in \{1,4,7,10,13,16,19,22,25\}$   $r = 2: x \in \{2,5,8,11,14,17,20,23,26\}$

$3 \mid x$   $X_{3,1} = \{0,1,2\}$              $r = 1: x \in \{3,12,21\}$                $r = 2: x \in \{6,15,24\}$

$3^2 \mid x$ $X_{3,0} = \{0\}$                  $r = 1: x \in \{9\}$                      $r = 2: x \in \{18\}$

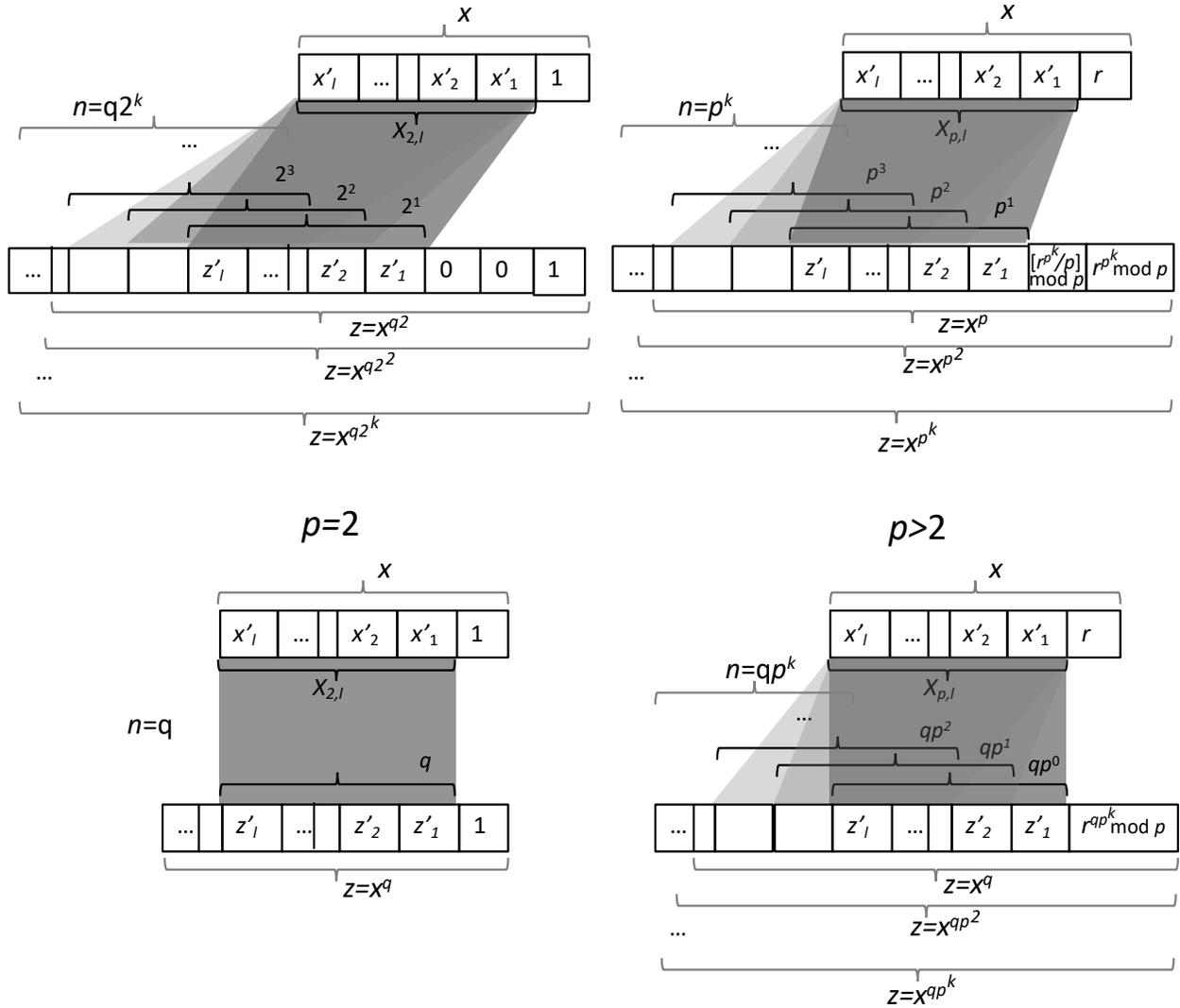

*Figure* 1. **The bijective mappings between consecutive runs of digits in the base $p$ representation of $x$ and $z = x^n$.** On the left side are the mappings base $p = 2$ and on the right are the mappings for a prime base $p > 2$. The lower left pane shows the simplest form of the mapping with coding shift 1 for $n = q$ with $q$ relatively prime to the base 2. There is a single non-zero residue class l mod 2 and for any integer $l > 0$ the next $l$ digits of $x$ take on all values base 2. The coding shift is 1 and coding region of the power $z = x^q$, depicted by the curly braces and shaded mapping, contains a permutation of numbers $0, 1, \ldots, 2^l - 1$ of the set $X_{2,l}$. When the power $n$ contains a power $k > 1$ of the base 2, the top left pane depicts the bijective mapping with coding shift $2 + k$ for $k = 1, 2, \ldots$. On the right are the corresponding cases for prime bases $p > 2$ for each non-zero residue class $r$ mod $p$. For clarity in these two diagrams, we include the relatively prime factor $q$ only in the lower pane and in the upper pane restrict the exponent $k$ of $p$ in the power $n = p^k$ to $k > 0$. The lighter shaded mappings are for successively higher exponents $k$, for which the coding regions in $z$ shift to the left.

For any power $n$, the extended PCT bijectively maps each of these six sets to equal size sets of numbers $z = x^n$. For example, for $n = 3$, and, for the first two sets where $3 \nmid x$, corresponding to the (non-extended) PCT, the two-digit coding region starts at the third digit ($\alpha = 2$) and the induced bijective mappings or permutations of $X_{3,2}$ for $r$ =1,2 are {0,7,2,3,1,5,6,4,8} and {0,4,2,3,7,5,6,1,8}. Note that the two lowest digits of the powers $z$ are not all distinct in this case with $z$ mod 9 being either 1 or 8 if $3 \nmid x$.

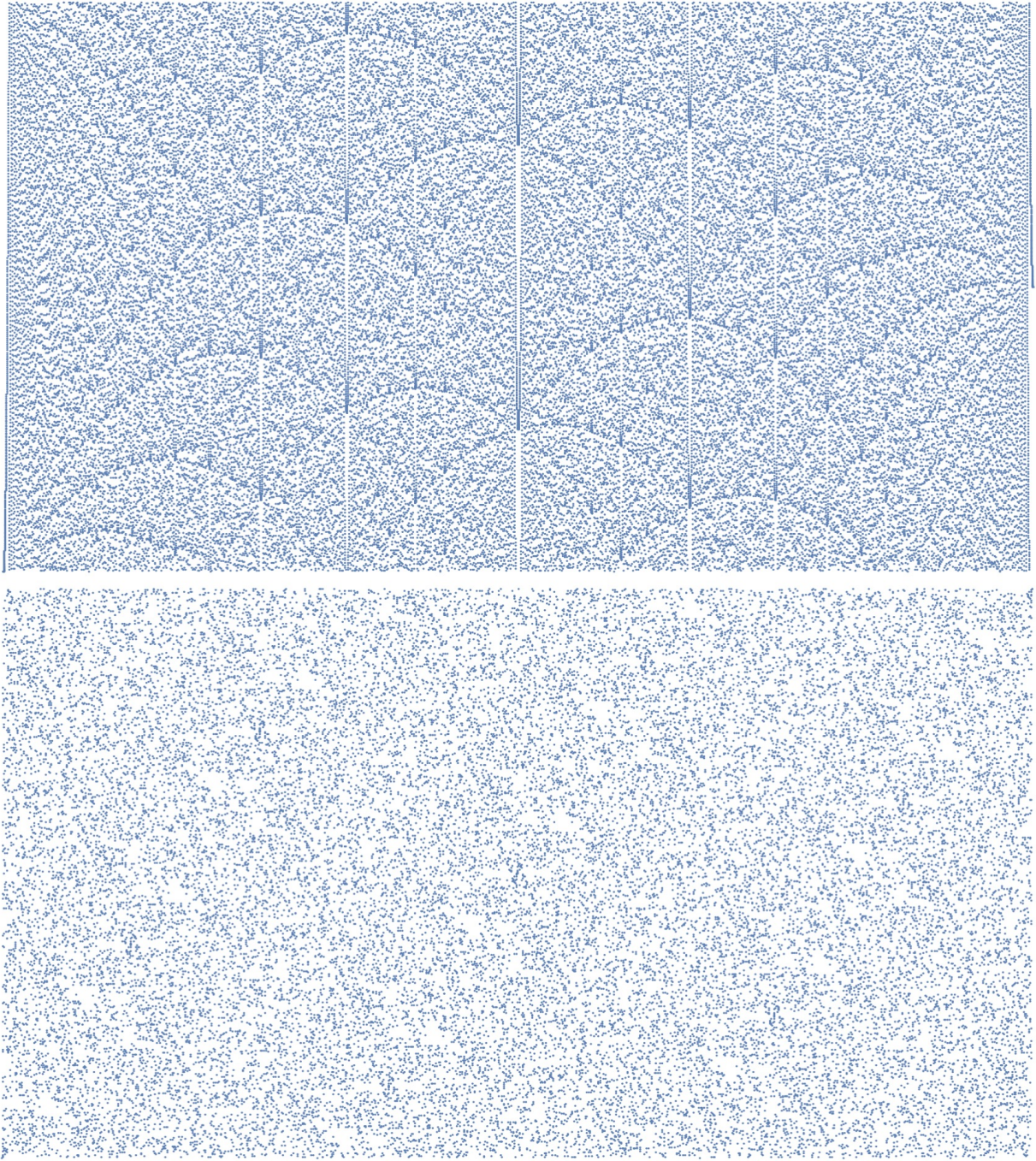

*Figure* **2. The bijective power coding mappings base $p$ = 2, power $n$ = 2, 3, radix $l$ = 16, 15 respectively.** The horizontal and vertical axes take on the values $x'$ and $z'$ respectively in $X_{2,l}$ in the range $0..\,2^l - 1$ and the points plot the mapping $C_1(x')$ for the single non-zero residue class 1 base 2. There are no collisions between points.

The mappings defined by the PCT can be depicted graphically, and we complete this section with some examples. In fig. 2 we plot two relatively large mappings for $p = 2$, with powers $n = 2,3$ showing how the power coding mapping maintains strict bijectivity for large sets $X_{2,l}$. This is achieved both when the mapping is highly structured as in the upper image ($n = 2$) or when it is seemingly random as in the lower image ($n = 3$).

## 4  Divisibility of Binomials involving Powers

We first prove a result about the divisibility of binomial coefficients involving prime powers:

<u>Lemma 1</u>

*For any prime p and positive integers j and k, if $p^h$ is the exact power of p dividing j, then $p^{k-h}$ is the exact power of p dividing the binomial coefficient $\binom{p^k}{j}$.*

<u>Proof</u>

We first observe that the integer binomial coefficients for $0 < j < p^k$ satisfy

$$\binom{p^k}{j} = \frac{p^k(p^k-1)\dots(p^k-j+1)}{j!} = \frac{p^k}{p^k-j}\binom{p^k-1}{j}$$

which means that if $j$ is not a multiple of $p$, i.e. $h = 0$, then $p^k - j$ is relatively prime to $p$ and hence to $p^k$, so that $p^k - j$ divides $\binom{p^k-1}{j}$ and thus more importantly $p^k$ divides $\binom{p^k}{j}$. On the other hand, if $j = g\,p^h$ with $g$ relatively prime to $p$ and $0 < h < k$, i.e. $j$ is a multiple of $p^h$, then we must prove that $p^{k-h}$ divides $\binom{p^k}{j}$.

Consider first for simplicity $h = 1$, meaning $j = gp$ and $p \nmid g$. Then $p^g$ divides $j!$, since the factorial $j! = \prod_{i=1}^{j} i$ in the denominator of $\binom{p^k}{j}$ increases by $p$ only when the index $i$ in the product is a multiple of $p$. But $(p^k - 1)\dots(p^k - j + 1)$ is divisible only by $p^{g-1}$ so that $p^{k+(g-1)-g} = p^{k-1}$ divides $\binom{p^k}{j}$. For $h > 1$, note that, incrementing $j$ from $p^h - 1$ to $p^h$, a new factor of $p^h$ appears in the denominator without compensation in the numerator (since $p^k - p^h + 1$ is relatively prime to $p$), so for $j = p^h$, $p^{k-h}$ divides $\binom{p^k}{j}$. At the next value $j = p^h + 1$, the last term in the numerator product is $p^k - j + 1 = p^k - p^h$ which returns to the binomial coefficient $\binom{p^k}{j}$ the factor of $p^h$ lost at the previous value $j = p^h$. We have established that $p^{k-h}$ divides $\binom{p^k}{j}$ for all $0 \leq h < k$.

Alternative proof *via* Kummer's Theorem (1852):

Kummer's theorem [5] gives an expression for the exact power of $p$ dividing the binomial coefficient $\binom{n}{m}$: it is the number of carries when adding $n - m$ and $m$ in base $p$. If $n = p^k$ with $k > 0$, then the base $p$ representation of $n = p^k$ is 1 followed by $k$ zeros. If $h$ is the exact power of $p$ dividing $m = j$, then the base $p$ representation of $j = qp^h$ is $(q_{k-h-1}, q_{k-h-2}, \ldots q_1, q_0, 0, \ldots, 0)$ with $h$ trailing zeros and $q_0 \neq 0$. The other base $p$ digits $q_i$ of $q$ are by definition in the range $0 \leq q_i < p$. So the base $p$ representation of $n - j$ is $(\check{q}_{k-h-1}, \check{q}_{k-h-2}, \ldots, \check{q}_1, \bar{q}_0, 0, \ldots, 0)$ with $\bar{q}_0 = p - q_0$ and $\check{q}_i = p - (q_i + 1)$. Note that the leading digits in this $k$ digit representation may be zero if $j < p^{k-1}$. Adding $n - j$ and $j$, the first carry occurs with $\bar{q}_0 + q_0 = p$ at index $h$ and then carries continue for every larger digit index since $\check{q}_i + q_i + 1 = p$ i.e. the number of carries is $k - h$. This proves the result. [Note that the carry count here is larger than the number of base $p$ digits where $n_i < m_i$ [6], since the latter does not include carry propagation for equal digits.]

## 5   Proof of the Power Coding Theorem

Our method of proof will be, for the general case of the power $n = q\,p^k$ with $q$ relatively prime to $p$,

(i)   to reduce the general case, as a composition of mappings, to the two separate cases (ii) and (iii):

(ii)  to prove the result for $n = q$ relatively prime to $p$

(iii) to prove the result for $n = p^k$, $p \geq 2$, $k \geq 1$. Because of the larger coding shift $\alpha > 1$, it is more convenient to prove the result directly for arbitrary $k$, than to try to employ an analogous composition of the mappings $n = p$ as in case (i).

For $n = 1$, the result is trivial, with the mappings being simply the identity. We henceforth assume $n \geq 2$.

**(i)**        $n = q\, p^k, p \geq 2, q$ **relatively prime to** $p$: **reduction to cases (ii) and (iii)**

This result follows from the fact that $x^{qp^k} = (x^q)^{p^k} = y^{p^k}$ where $y = x^q$. If case (ii) is true for $x^q$, writing $x = px' + r$ with $x' \in X_{l,p}$, there is, for each non-zero residue $r$ base $p$, a bijective mapping $f:X_{l,p} \to X_{l,p}$ mapping $x'$ to $y'$ with coding shift $\alpha = 1$ where $y = x^q \bmod p^{l+1} = py' + r'$ where $r' = r^q \bmod p$. Since the shift is one in this case, and since $r' \neq 0$ because $q$ and $p$ are relatively prime and $r \neq 0$, the set of numbers $y$ resulting from the bijective mapping is precisely of the form $py' + r'$, with $y' \in X_{l,p}$ and $r' \neq 0$, required for renewed application of the theorem without any overall shift in the base $p$ representation. If case (iii) is true, the bijective mapping $g:X_{l,p} \to X_{l,p}$ based on $y^{p^k}$, as defined and proved as the bijective mapping in (iii) below, with shift $\alpha_g$, then completes the composite mapping $C_r = gf:X_{l,p} \to X_{l,p}$ with total coding shift $\alpha_g$. So the proof of this case follows directly from cases (ii) and (iii) below.

**(ii)**        $n = q, q$ **relatively prime to** $p$.

For $l > 0$ in the theorem, let $N = p^{l+1}$, then the natural numbers $x = x'p + r$, with $x' \in X_{l,p}$ (and $0 < r < p$), are relatively prime to $N$ and hence form a group under multiplication modulo $N$ (usually referred to as the group of units modulo $N$, $\mathbb{Z}_N^*$). The number of such units, i.e. the order of the group, is the Euler Totient function [2] $\phi(N) = \phi(p^{l+1}) = p^l(p-1)$. This group has a ($p$-Sylow) subgroup, $G_{p,l}$, of order $p^l$, consisting of the elements of $\mathbb{Z}_N^*$ with $r = 1$.

As pointed out to us by H. Braden (personal communication), Gauss showed [7] that, in our notation, unless the prime base $p = 2$ and $l > 1$, $\mathbb{Z}_N^*$ and hence $G_{p,l}$ is also cyclic, $G_{p,l} = C_{p^l}$. The restriction on $p$ and $l$ makes it more convenient to prove the result for all cases without recourse to this simplifying feature. This also makes the proof more direct. Since $q$ is relatively prime to $p$, and hence to $p^{l+1}$, there exists an integer $s > 0$ satisfying $qs \equiv 1 \pmod{N}$, i.e. for some integer $t > 0$, $qs = tN + 1$. Then $x(z) = z^s \bmod N$ demonstrates an inverse function for $z(x) = x^q \bmod N$. To see this, we note that

$$z(x)^s \equiv x^{qs} \equiv x^{tN+1} \equiv x\left(x^{p^l}\right)^{tp} \equiv x \pmod{N}$$

since $x^{p^l} \equiv 1 \pmod{p^{l+1}}$ by Euler's theorem for the group $G_{p,l}$. Since the mapping $z(x)$ is also onto, we have proven that it is a permutation of $G_{p,l}$. It follows that $\alpha(q,p) = 1$ and that the mapping $C_1(x') = \dfrac{(px'+1)^q - 1}{p} \bmod p^l$ is bijective. Note that the elements of $\mathbb{Z}_{p^{l+1}}^*$ are partitioned into the cosets of $G_{p,l}$ in $\mathbb{Z}_N^*$, namely $rG_{p,l}$ for $r > 0$, each of which is in 1:1 correspondence to $G_{p,l}$, and consists of the $p^l$ distinct elements with remainder $r \bmod p$. Furthermore, we could equally have chosen any representatives of these $p-1$ cosets to form this partition of $\mathbb{Z}_N^*$ from $G_{p,l}$. In particular, the cosets $r^q G_{p,l}$ with multiplication in $\mathbb{Z}_N^*$ would also do, provided that the $p-1$ elements $r^q$ mod $p$ are distinct. However this is only true if $p-1$ is relatively prime to $q$. With this assumption, the cosets partition $\mathbb{Z}_N^*$ and our result for $G_{p,l}$ implies that the mapping $x^q$ mod $N$ is 1:1 on the whole of $\mathbb{Z}_{p^{l+1}}^*$. In general, the power coding mappings are only bijective separately on the individual residue cosets $rG_{p,l}$.

We have shown that the mapping of the form eqn(1.1) $C_r(x') = \frac{(px'+1)^q - (r^q \bmod p)}{p}$ mod $p^l$ is bijective and the result is proven for $n = q$ relatively prime to $p$.

(iii)     $n = p^k$, $p \geq 2$, $k \geq 1$

As above, let $x = px' + r$, with $x' \in X_{l,p}$ and let $N = p^{l+k+\delta_{p,2}+1}$. Consider the mappings

$$z(x') = x^{p^k} \bmod N = (px' + r)^{p^k} \bmod N$$

$$y(x') = p^{k+\delta_{p,2}+1} x' + \left(r^{p^k} \bmod p^{k+\delta(p,2)+1}\right) \bmod N.$$

The two terms in the mapping $y(x)$ separate the first $k + \delta(p, 2) + 1$ digits base p, constant as a function of $x'$, from the higher digits. The mapping $y$ is 1:1, with inverse defined on its range

$$x(y) = p \lfloor y/p^{k+\delta(p,2)+1} \rfloor + r$$

since the lower digits are constant and the upper digits are simply a shifted copy of $x'$, so that we can consider the induced mapping $y \to z$ defined uniquely on the range of $y$. We claim that $z(y)$ is a permutation of the range of $y(x')$.

Using lemma 1, we can now show that $z(x)$ maps $x$ onto the range of $y$. For

$$z(x') = (px' + r)^{p^k} = \sum_{j=0}^{p^k} \binom{p^k}{j} p^j x'^j r^{p^k-j}$$

$$= r^{p^k} + p^k p\, x'\, r^{p^k-1} + \frac{p^k(p^k-1)}{2} p^2 x'^2\, r^{p^k-2} + \sum_{j=3}^{p^k} \binom{p^k}{j} p^j x'^j r^{p^k-j}$$

$$= \begin{cases} 1 + 2^{k+1}(x' + (2^k-1)x'^2) + 2^{k+2} \sum_{h=0}^{k} \sum_{j=3, j=g2^h, 2\nmid g}^{2^k} \frac{\binom{2^k}{j}}{2^{k-h}} 2^{j-2-h} x'^j, & p = 2 \\[2ex] r^{p^k} + p^{k+1} x'\, r^{p^k-2}\left(r + p\frac{(p^k-1)}{2} x'\right) + p^{k+1} \sum_{h=0}^{k} \sum_{j=3, j=gp^h, p\nmid g}^{p^k} \frac{\binom{p^k}{j}}{p^{k-h}} p^{j-1-h} x'^j r^{p^k-j}, & p > 2 \end{cases}$$

$$\equiv r^{p^k} \pmod{p^{k+\delta_{p,2}+1}}$$

In the third line of this derivation, we have separated the cases $p = 2$ and $p > 2$, and in the sum we have regrouped together all terms with the same exact power $p^h$ of $p$ in $j$. We have grouped terms so that all exponents appearing in these equations are non negative and all fractional expressions are integers. Especially, since $j \geq p^h > h$, $j - 1 - h \geq 0$, and for $p = 2$ and $j \geq 3$ we have $j - 2 - h \geq 0$. The binomial divisibility in Lemma 1 establishes that $p^{k-h} | \binom{p^k}{j}$ for $j = g\, p^h$ where $p \nmid g$, and this completes the proof that all the terms in the double sum are integers. Furthermore, when $p > 2$, $p^k - 1$ is even and so the second term is an integer and multiple of $p^{k+1}$. For the case $p = 2$, note that $x'$ and $x'^2$ have the same parity and so the terms in parentheses are both even or both odd, so that their sum is even, raising the power of 2 in the term to $2^{k+2}$ and allowing the last line to be deduced using the Kronecker delta. So the mapping $z$, is onto the range of $y$.

It remains to prove that the mapping $z$ is 1:1. For $l \geq 1$, let $l' = l + k + \delta_{p,2} + 1$, and note that $\phi(p^{l'}) = p^{l'-1}(p-1) > p$ and that by the Fermat-Euler theorem [2], since $p \nmid x$, $x^{\phi(p^{l'})} = 1 \bmod p^{l'}$. We can construct an inverse mapping to $z$ using multiplication mod $p^{l'}$ by $u(x) = x^{\phi(p^{l'})-p^k+1} \bmod p^{l'}$, since $x^{p^k} x^{\phi(p^{l'})-p^k+1} \bmod p^{l'} = xx^{\phi(p^{l'})} \bmod p^{l'} = x \bmod p^{l'}$. The existence of an inverse mapping establishes that it is also 1:1. This arguments holds for all $l \geq 1$, so we have completed the proof of case (iii).

This completes the proof of the theorem.

## 6   Discussion and Conclusions

The power coding mapping is obviously related, through the use of the power mapping, to the Frobenius automorphism for a finite field $F$ of characteristic $p \neq 0$, given by the map $\phi : F \to F$ which maps each element $u$ of $F$ to $u^p$. However, we are not dealing here with finite fields, but particular subsets of numbers embedded within the base p encodings of numbers and their powers. Moreover, the power coding mapping is defined for arbitrary $n$, independently of the choice of prime base $p$.

It has not escaped the authors, that the 1:1 recoding of entire sets of numbers using powers (which already has applications in cryptography) potentially may have a bearing on simple proofs of Fermat's Last Theorem. This is however being investigated separately.

**Acknowledgements**

The theorem was first stated as a numerical conjecture by the authors in April 2001 during a research visit of P.R.W., funded by the Humboldt Foundation, to the Research Group BioMIP at the German National Research Center for Information Technology (GMD). Initial symbolic algebra verification and graphical display of the bijective mappings was carried out using *Mathematica™*. The authors wish to thank Andreas Dress for reading and commenting on the manuscript and Harry Braden for helpful comments and proof reading the equations. The first author wishes to dedicate the theorem to Elizabeth Esther McCaskill, who taught him about number relationships in colour and lengths at an early age.

Appendix 1 Alternative proof of PCT for $n = p = 2$

In this special case, let $x = 2x' + 1$, with $x' \in X_{l,2}$ and let $N=2^{l+3}$. Consider the mappings:

$$z(x') = x^2 \bmod N = (2x' + 1)^2 \bmod N \text{ and } y(x') = 4x - 3 = 8x' + 1, \text{ for } x' \in X_{l,2}$$

noting that $4x - 3 < N$. The linear mapping $y(x')$ is 1:1, so that we can consider the induced mapping $y \to z$ defined uniquely on the range of $y$, since $y(x')$ being 1:1 is invertible on the range of $y$ (in fact $x'(y) = \lfloor y/8 \rfloor$). We claim that $z(y)$ is a permutation of the range of $y(x)$. Firstly, $z(x)$ maps $x$ onto the range of $y$, for

$$z(x') = x^2 \bmod N$$
$$= (2x' + 1)^2 \bmod N$$
$$= (4x'^2 + 4x' + 1) \bmod N$$
$$= (8x'(x' + 1)/2 + 1) \bmod N$$

since $x'(x' + 1)$ is even for all $x'$. So $z(y)$ is also onto for the range of $y$. Secondly, assume that the mapping $z(y)$ is not 1:1. We demonstrate that this leads to a contradiction. For let $x_1 = 2x_1' + 1$ and $x_2 = 2x_2' + 1$ be two different positive integers with $x_1', x_2' \in X_{l,2}$ and $x_1' > x_2'$ satisfying

$$z(x_1) = z(x_2)$$
$$x_1^2 \equiv x_2^2 \pmod{N}$$
$$4x_1'(x_1' + 1) \equiv 4x_2'(x_2' + 1) \pmod{N}$$

so that by assumption

$$z(x_1) - z(x_2) \equiv 4\big(x_1'(x_1' + 1) - x_2'(x_2' + 1)\big) \pmod{N}$$
$$\equiv 4\big((x_1'^2 - x_2'^2) + (x_1' - x_2')\big) \pmod{N}$$
$$\equiv 4\big((x_1' - x_2')(x_1' + x_2' + 1)\big) \pmod{N}$$
$$\equiv 0 \pmod{N}$$

If $x_1'$ and $x_2'$ are both even or both odd, then the second factor in line three is odd and relatively prime to $N = 2^{l+3}$ so that it may be cancelled. Then $x_1' = x_2'$, since 4 times their difference is less than $N$. This contradicts our assumption. If, on the other hand only one of $x_1'$ and $x_2'$ is even, then $0 = 4(x_1' + x_2' + 1) \bmod N$, since their difference is odd and may be cancelled. But $4(x_1' + x_2' + 1) < N$ so that $x_1' + x_2' + 1 = 0$, which cannot be satisfied by $x_1', x_2' \in X_{l,2}$, and so we also have a contradiction. Hence the mapping $z(x)$ is 1:1, and the mapping $z(y)$ is 1:1 and onto for the range of $y(x)$. It follows that $\alpha(2,2) = 3$ and that the mapping $C_1(x') = \frac{(2x'+1)^2 - 1}{2^3} \bmod 2^l = \frac{1}{2}x'(x' + 1) \bmod 2^l$ is bijective. A closely related mapping was proposed for random number generation by R.R. Coveyou [8].

## Appendix 2: Alternative proof of PCT for $n = 2^k$, $p = 2$

Let $x = 2x' + 1$, with $x' \in X_{l,2}$ but now let $N = 2^{l+k+2}$. Consider the pair of mappings $z_k(x') = x^{2^k} \bmod N = (2x' + 1)^{2^k} \bmod N$ and $y_k(x') = 2^{k+2}x' + 1$. The linear mapping $y_k(x)$ is 1:1, so that we can consider the induced mapping $y \to z$ defined uniquely on the range of $y$. We show that $z_k(y)$ is a permutation of the range of $y_k(x)$, proceeding by induction on $k$. By (ii) the result is true for $k = 1$. We assume it is true for $k$, and prove the result for $k + 1$. We must first prove that $z_{k+1}(y)$ is onto, for the range of $y_{k+1}(x)$. This follows directly from the result for $k$, since

$$(2x' + 1)^{2^{k+1}} - 1 = \left((2x' + 1)^{2^k} - 1\right)\left((2x' + 1)^{2^k} + 1\right)$$

and the first factor is divisible by $2^{k+2}$ (since the result is assumed true for $k$) and the second factor is divisible by 2, being the sum of two odd numbers. Thus $z_{k+1}(x) - 1$ is divisible by $2^{k+3}$, establishing the inductive result. Next we assume that $z_k(x)$ is 1:1 and must prove that this is true also for $z_{k+1}(x)$. We note that

$$((2x' + 1)^{2^k}) \bmod 2^{l+k+3} = a(x', k, l)2^{l+k+2} + z_k(x)$$

where $a(x', k, l)$ is either 0 or 1. Looking again at the mapping $z_{k+1}(x)$,

$$(2x' + 1)^{2^{k+1}} \equiv z_k(x)^2 \pmod{2^{l+k+3}}$$

All terms involving $a(x', k, l)$ contain powers of 2 higher than $2^{l+k+2}$ and hence are zero (mod $2^{l+k+3}$). Now this is the composition of two 1:1 mappings, $z_k(x)$ and $z_1(x)$ (the latter for a different value of $l$) by the inductive assumption and hence is also 1:1. This completes the proof for this case.

Appendix 3: Alternative proof of PCT for $n = p$, $p > 2$, $r = 1$

Let $x = px' + 1$, with $x' \in X_{l,p}$ and let $N = p^{l+2}$. We consider the mappings $z(x') = x^p \bmod N = (px' + 1)^p \bmod N$ and $y(x') = p(x - 1) + 1 = (p^2 x' + 1) \bmod N$. The linear mapping $y(x')$ is 1:1, so that we can consider the induced mapping $y \to z$ defined uniquely on the range of $y$. In fact $y(x)$, being 1:1, is invertible on the range of $y$ (in fact $x(y) = \lfloor y/p^2 \rfloor$). We claim that $z(y)$ is a permutation of the range of $y(x)$. Firstly, $z(x')$ maps $x'$ onto the range of $y$, for

$$z(x') = x^p \bmod N$$
$$= (px' + 1)^p \bmod N$$
$$= \sum_{j=0}^{p} \binom{p}{j} p^j x'^j \bmod N$$
$$= \left(1 + p^2 \sum_{j=1}^{p} \binom{p}{j} p^{j-2} x'^j\right) \bmod N$$

since $\binom{p}{j}$ is divisible by $p$ for $1 \leq j \leq 2$ (actually up to $p - 1$). So the mapping is onto the range of $y(x)$.

To establish that the mapping is also 1:1 for $r = 1$, we proceed inductively on $l$ for $l > 0$. Firstly, for $l = 1$, we have $x = px' + 1$ with $0 \leq x' < p$. But $x'^p = x' \bmod p$, so that $x^p = (p^2 x' + 1) \bmod p^3$ which is simply $y(x')$ and clearly a 1:1 mapping of $x'$. Thus the result is true for $l = 1$. Now assume the mapping is 1:1 for $l$, we must prove it is 1:1 for $l+1$. For $x \in X_{l+2,p}$, we may write $x = x_{l+1} p^{l+1} + x'_l p + 1$ where $x'_l \in X_{l,p}$ and $0 \leq x_{l+1} < p$, so that

$$z_{l+1}(x) = x^p \bmod p^{l+3}$$
$$= \left(x_{l+1} p^{l+1} + (x'_l p + 1)\right)^p \bmod p^{l+3}$$
$$= ((x'_l p + 1)^p + p(x'_l p + 1)^{p-1} x_{l+1} p^{l+1}) \bmod p^{l+3}$$
$$= (z_l(x'_l) + x_{l+1} p^{l+2}) \bmod p^{l+3}$$

where in the third line higher powers of $x_{l+1} p^{l+1}$ in the binomial expansion are zero mod $p^{l+3}$ since $2(l + 1) = 2l + 2 > l + 3$ for $l > 1$. Similarly, in the last line, only the first term 1 in the binomial expansion of $(x'_l p + 1)^{p-1}$ is needed. Since the first term mapping $z_l(x'_l)$ is 1:1 and $z_l(x'_l) < p^{l+2}$ by our inductive hypothesis so that $z_l(x'_l) + x_{l+1} p^{l+2} < p^{l+3}$ and hence the sum of two 1:1 mappings on $\mathbb{Z}$ and hence 1:1. This completes the inductive proof for this case with $r = 1$.